\documentclass[11pt,a4paper]{article}
\usepackage{amsfonts,amsthm}
\usepackage{fullpage}

\newcommand{\N}{\mathbb{N}}
\newcommand{\erre}{\mathbb{R}}

\newcommand{\ds}{\displaystyle}

\renewcommand{\epsilon}{\varepsilon}

\newcommand{\cp}[2]{\left\langle#1,#2\right\rangle}

\newtheorem{prop}{Proposition}%[section]

\newtheorem{lemma}[prop]{Lemma}
\theoremstyle{remark}
\newtheorem{rmk}[prop]{Remark}

\title{Optimal Distributed Dynamic Advertising}

\author{Carlo Marinelli\thanks{
Institut f\"ur Angewandte Mathematik, Universit\"at Bonn,
Wegelerstr. 6, D-53115 Bonn, Germany. E-mail
\texttt{marinelli@wiener.iam.uni-bonn.de}. Corresponding author.}
\ and Sergei Savin\thanks{
Graduate School of Business, Columbia University,
3022 Broadway, New York, NY 10027, USA}}

\date{\normalsize November 2006}

\begin{document}
\maketitle

\begin{abstract}
We propose a novel approach to modeling advertising dynamics for a
firm operating over distributed market domain based on controlled partial
differential equations of diffusion type. Using our model, we consider 
a general type of finite-horizon profit maximization problem in a monopoly  
setting. By reformulating this profit maximization problem as an optimal control    
problem in infinite dimensions, we derive sufficient conditions for the     
existence of its optimal solutions under general profit functions, as well  
as state and control constraints, and provide general characterization of   
the optimal solutions. Sharper, feedback-form, characterizations of the     
optimal solutions are obtained for two variants of the general problem.

\medskip\par 


\medskip\par

\end{abstract}

\section{Introduction}
In the last two decades, dynamic optimal control models have become a
mainstay in the marketing literature focused on advertising. While the
empirical marketing research pays increasingly growing attention to
the multiple-market aspects of the advertising strategy (\cite{bronn},
\cite{dube}), the modeling support for the studies of advertising
policies distributed over multiple geographic regions has been lacking
(see \cite{sethi} for the most comprehensive review of the related
literature).  In the present paper, we propose a simple model which
deals with both the dynamic as well as the spatial effects of
advertising.

The classical paper by Nerlove and Arrow \cite{NA} introduced the following
model for the dynamics of goodwill stock under the influence of advertising
for a monopolist in a single market-single product environment: 
\begin{equation}
{\frac{dx(t)}{dt}}=u(t)-\rho x(t),  \label{eq:NA}
\end{equation}%
where $x(t)$ is the goodwill level at time $t\geq 0$, $u$ is the rate
of advertising (in monetary terms), and $\rho>0$ is a constant factor
describing the deterioration of the goodwill in the absence of
advertising.  Our analysis extends the Nerlove-Arrow (NA hereafter)
model by considering spatially distributed advertising and by modeling
geographic fluctuations in the goodwill stock. The treatment of
problems of optimal advertising in the space-time setting is not new:
to the best of our knowledge it was considered for the first and only
time in \cite{SSD} as a special case of a very general class of
distributed-parameter optimal advertising models, where the parameter
space is modeled as subset of a measurable space. The aim of
\cite{SSD} was to establish a general, abstract framework and
essentially proved the well-posedness of the problem and the existence
of optimal strategies under suitable assumptions. Our purpose here is
different, as we focus instead on a simple NA model with a diffusive
component, for which optimal strategies can be better characterized,
sometimes even explicitly.

In the first part of this work we address a modeling question focused
on describing the dynamics of goodwill which depends on both space and
time coordinates. In particular, we propose to model the dynamics of
goodwill through the following controlled partial differential
equation (PDE)
\begin{equation}
\label{eq:PDE}
{\frac{\partial x}{\partial t}}(t,\xi)=-\rho x(t,\xi)
+ \Delta_\xi x(t,\xi) + b(\xi)u(t,\xi)
\end{equation}
(with appropriate initial and boundary conditions -- see
(\ref{eq:pde}) below), where $x:[0,T]\times \Xi \to \erre$, $\Xi
\subset \erre^2$ is the goodwill ``density'', $u(\cdot,\cdot)$ is the
rate of advertising effort (expressed in Gross Rating Points, or
GRPs), and $b(\cdot)$ is the coefficient of advertising effectiveness.
The second term on the right-hand side of (\ref{eq:PDE}) is introduced
to capture the effect of the spatial diffusion of goodwill (more
detailed rationale for the use of this term is described in the next
section) and reflects the main modeling difference between our
analysis and that of Nerlove and Arrow \cite{NA}. The last term on the
right-hand side of (\ref{eq:PDE}) represents another new feature of
our model: we describe advertising effort using GRP parameter, which
recent trend in advertising literature\ (\cite{vilca}, \cite{dube})
defines as more appropriate than the traditional expenditure rate.
Note that (\ref{eq:PDE}) reduces to the NA model if $x$ does not
depend on the spatial coordinate $\xi$ and if $b(\xi)\equiv 1$.
Moreover, we shall see that the NA model is also recovered by
averaging in space our model.

In the second part, we formulate and solve the space--time analogs of
some of the optimal control problems studied earlier for the cases
where $x$ depends only on time. Most of the results will follow using
the tools of optimal control in infinite dimensional spaces.
We should clearly state at this point that the emphasis in the paper
is not on deriving new results in infinite dimensional optimal
control, but rather on showing that models and techniques of analysis
and control in Hilbert spaces could prove useful to extend existing
advertising models to account for geographic fluctuations. Moreover,
we devote some effort to show that in some simple, but still
interesting situations, optimal policies can be obtained in closed
form, and are hence amenable to practical implementation. Finally, we
discuss how the diffusive term influences the optimal strategy by a
direct comparison with the standard NA model.

In particular, we consider the following types of problems: in the
most general case our aim is to maximize a monopolistic firm's finite
horizon profit expressed as a combination of a functional of the
terminal goodwill stock and the cumulative cost of advertising. We
approach this problem in the abstract setting of control of
infinite-dimensional systems via Pontryagin's maximum principle: under
mild assumptions on the structure of the problem, one can obtain quite
general, although abstract, results on the existence and
characterization of optimal policies. For some specific choices of
cost and reward functions, optimal policies can be expressed in closed
form. A dynamic programming approach allows one to characterize the
value function as solution of an Hamilton-Jacobi equation, and to
express optimal strategies in feedback form. In the special case of
quadratic cost and reward functions, a linear-quadratic (LQ) regulator
approach yields explicit expressions for optimal strategies in
feedback form.  The LQ approach is also used to solve a related
problem, i.e. the minimization of a weighted sum of the distance of
$x(T,\cdot)$ from a target goodwill level and the total cost of
reaching the target.

Our analysis follows a classical scheme: first we establish
that the controlled PDE (\ref{eq:PDE}) can be written as an abstract
linear control system of the type
\begin{equation}
\frac{dx}{dt}=Ax(t)+Bu(t)  \label{eq:ODE}
\end{equation}%
in a suitable Hilbert space of functions $X$. Then we show that the
control problem for the PDE (\ref{eq:PDE}) is equivalent to a control
problem on $X$ for the abstract differential equation (\ref{eq:ODE}).
In the general case, the optimal advertising problem is solved
(although only in an abstract way) using the weak maximum principle in
infinite dimensions, the theory of which can be found, e.g., in
\cite{barbu}. The simpler case of quadratic utility and cost
functions is reduced to the study of a Riccati equation by an
application of the dynamic programming principle. The solution of such
an equation yields, in particular, a feedback-type optimal policy for
the original problem. For the infinite-dimensional LQ problem we refer
to \cite{rep1}, \cite{rep2}, \cite{CZ} for the standard case of
positive definite costs (see also \cite{LT1}), and to \cite{LY} for
the general case with indefinite costs.

The rest of the paper is organized as follows: section \ref{sec:model}
formally defines and motivates the model for the dynamics of goodwill
as a function of space and time, as well as the associated optimal
advertising problems. The model as well as the optimization problems
are recast in the framework of control systems in Hilbert spaces in
section \ref{sec:reform}. Sections 4 to 7 deal with the optimal
control problems mentioned above, and section 8 concludes.

\section{Model for spatial advertising dynamics}
\label{sec:model}
Let us first introduce some notation used throughout the paper. Let
$X$ be a real separable Hilbert space with norm $|\cdot|$ and inner
product $\langle\cdot,\cdot\rangle$. We denote by $L_2([0,T],X)$
the Hilbert space of functions $f:[0,T]\to X$ with
finite norm $\Vert f\Vert $,
$$
\Vert f\Vert ^{2}=\int_{0}^{T}|f(t)|^{2}\,dt<\infty .
$$
When there will be no possibility of confusion, we shall still denote
by $\cp{\cdot}{\cdot}$ the inner product of spaces like
$L_2([0,T],X)$. The space of linear bounded operators between the
Hilbert spaces $X$ and $Y$ will be denoted as $\mathcal{L}(X,Y)$, or
simply as $\mathcal{L}(X)$ when $X=Y$.

Given an operator $A:D(A)\subset X\to X$, we shall say that
$A$ is {\em uniformly positive definite} (denoted $A>\!\!\!>0$) if there
exists $\varepsilon>0$ such that $A-\varepsilon I\geq 0$, i.e.
$\cp{Ax}{x} \geq \varepsilon |x|^2$ for all $x\in D(A)$.
The adjoint of an operator $A$ is denoted by $A^*$.

We shall need some standard functional spaces: for a bounded open set
$\Omega \subset \erre^d$, $H^k(\Omega )$ is the Sobolev space of
functions with (generalized) derivatives up of order $k$ in
$L_2(\Omega)$, and $H_{0}^{1}(\Omega)$ is the closure of
$C_{0}^{\infty }(\Omega )$ in the topology of $H^{1}(\Omega)$. Here
$C_0^\infty(\Omega)$ stands for the space of infinite differentiable
functions with compact support on $\Omega$.

The following elementary notions from convex analysis will be used: given a
Hilbert space $X$ and convex lower semicontinuous function $f:X\to 
]-\infty ,+\infty ]$, the mapping $\partial f:X\to X$ defined by   
$$
\partial f(x)=\{w\in X:\langle x-y,w\rangle \geq f(x)-f(y),\forall y\in X\}
$$
is called the subdifferential of $f$. Moreover, the conjugate (or
Fenchel-Legendre transform) of $f$ is the function $f^*:X\to
]-\infty ,+\infty ]$ defined by 
$$
f^{\ast }(p)=\sup_{x\in X}(\langle p,x\rangle -f(x)).
$$

\medskip

Let $\Xi$ be a bounded open set of $\erre^2$ with regular boundary
$\partial\Xi$. This set will be the model for the geographic region of
interest for the advertising campaign. Let $x:[0,T]\times\Xi\to\erre$
be the goodwill ``density'' of a given product, specified as a
function of time $t$ and location $\xi \in \Xi $. The model for the
controlled dynamics of $x$ will be given by the following equation:
\begin{equation}
\label{eq:pde}
\left\{ 
\begin{array}{ll}
\ds \frac{\partial x}{\partial t}(t,\xi)=(-\rho +\Delta_\xi)x(t,\xi )
+b(\xi)u(t,\xi), & (t,\xi )\in ]0,T]\times \Xi , \\[10pt] 
\ds \frac{\partial x}{\partial n}(t,\xi)=0, &
           (t,\xi)\in ]0,T]\times \partial\Xi , \\[8pt] 
x(0,\xi )=x_{0}(\xi ), & \xi \in \Xi ,%
\end{array}%
\right.  
\end{equation}%
where $T>0$ is a fixed time (which could be thought, for instance, as the
time of market introduction for the new product), $\rho >0$ is a natural 
deterioration factor of the product image in absence of advertising, $\Delta
_{\xi }=\frac{\partial ^{2}}{\partial \xi _{1}^{2}}+\frac{\partial ^{2}}{%  
  \partial \xi _{2}^{2}}$ is the Laplacian with respect to the spatial
variable $\xi $, $u:[0,T]\times \Xi \to \erre_{+}$ is the
advertising effort, which also depends on time and location, and
$b:\Xi \to \erre_{+}$ is a (bounded) factor of advertising
effectiveness. $x_{0}:\Xi \to \erre$ represent the initial
level of goodwill as a function of location in the region of interest
(for example, one may assume $x_{0}=0$ if at time $t=0$ the advertised
product is unknown).

The introduction of the diffusion term $\Delta_{\xi} x(t,\xi)$ is
motivated by an analogy with the classical diffusion equation, used to
model the evolution in time of the concentration of a substance in an
homogeneous medium. The idea is to regard a consumer population in a
bounded region $\Xi$ as the medium in which the goodwill is spreading.
The Neumann boundary condition in (\ref{eq:pde}) is equivalent to
assuming that there is no exchange of goodwill through the boundary of
the region $\Xi$, i.e. that the consumer population inside $\Xi$ does
not communicate with the outside world. This assumption is natural (at
least when considering large enough populations, such as a nation) and
allow us to recover the classical NA by averaging our model in space.
Indeed, integrating (\ref{eq:pde}) we get
$$
\int_\Xi \frac{\partial x}{\partial t}(t,\xi)\,d\xi = 
-\rho \int_\Xi x(t,\xi)\,d\xi + \int_\Xi \Delta x(t,\xi)\,d\xi
+ \int_\Xi b(\xi)u(t,\xi)\,d\xi.
$$
The total goodwill in the region $\Xi$ is given by
$\bar{x}(t):=\int_\Xi x(t,\xi)\,d\xi$. Set $\bar{u}(t)=\int_\Xi
b(\xi)u(t,\xi)\,d\xi$ (in the special case $b=1$ this would be the
total advertising effort).  Note that, by the Gauss-Stokes formula and
the Neumann boundary conditions in (\ref{eq:pde}), one has
$$
\int_\Xi \Delta x(t,\xi)\,d\xi =
\int_\Xi \Delta x(t,\xi)1(\xi)\,d\xi =
\int_{\partial\Xi} \nabla x(t,\xi)\nabla 1(\xi)\,d\xi = 0.
$$
Hence, by changing the order of differentiation and integration on
the left hand side,
% (and this can be done, e.g. by a simple application of
% Fatou's lemma),
we obtain
$$
\frac{d\bar{x}}{dt} = -\rho\bar{x} + \bar{u},
$$
which represents the classical NA dynamics for the total goodwill
$\bar{x}(\cdot)$.

\begin{rmk}
  (i) Equations like (\ref{eq:pde}) have also been successfully used
  to model the dynamics of populations, as well as the dynamics of
  infective diseases, in bounded regions in $\erre^2$ or $\erre^3$.
  Let us also recall that simple epidemics models (without spatial
  dimension) were the starting point for Bass' model of the diffusion of new
  goods, which still is a reference model in the marketing and
  management science literature. These considerations, together with
  the above ones, would hopefully give a plausible foundation to our
  model.

  (ii) The concept of ``goodwill density'' is not at all difficult to
  interpret nor too abstract. In fact, it can be (roughly) interpreted
  as the ratio of total goodwill over a region and the (2-dimensional
  Lebesgue) measure of the region, provided the region is small
  enough.

  (ii) In the case of a non-homogeneous consumer population the term
  $\Delta x(t,\xi)$ modeling goodwill diffusion would be replaced by a
  more general second-order elliptical operator, adding only some
  technical complications.
\end{rmk}

Below we will use the controlled dynamic model (\ref{eq:pde}) to
identify the optimal advertising strategy for a monopolistic firm.
Traditional advertising literature provides several choices in
selecting the objective to be maximized in such setting. In the
present work we adopt the perspective of a firm which is prepared to
launch a new product or service at a pre-determined time $T$ in the
future and which aims at maximize its total pre-launch profits.
Therefore, a rather general distributed advertising problem one would
like to solve can be described as follows:
\begin{itemize}
\item[(P)] maximize the functional 
\begin{equation}
J_{c}(u)=\int_{\Xi }\phi _{0}(x(T,\xi ))\,d\xi -\int_{0}^{T}\!\!\int_{\Xi
}h_{0}(u(t,\xi ))\,d\xi \,dt  \label{eq:general}
\end{equation}%
over all controls $u(t,\xi )\in [0,R]$ $(t,\xi)$-a.e., subject to
the dynamics (\ref{eq:pde}) and the additional constraint $x(t,\xi )\geq 0$
$(t,\xi)$-a.e.
\end{itemize}
Here $\phi_0$, $h_0:\erre\to\erre$ are such that the
integrals in (\ref{eq:general}) are finite, and can be thought as
utility (profit projection) of final goodwill and unit cost of
pre-launch advertising, respectively.

We shall see that if the cost function $h_0$ is quadratic, an
alternative approach based on Hamilton-Jacobi equations allows one to
obtain a feedback characterization of the optimal strategy. Quadratic
functions provide the simplest (strictly) convex cost functions and are
commonly used (see e.g. \cite{muller}).

\smallskip

We shall also consider two simpler variants of the general advertising
problem for which sharper characterizations of the optimal advertising
policy can be obtained. The first one, which we denote as (P1), is
formulated as follows:
\begin{itemize}
\item[(P1)] maximize the functional 
\begin{equation}
\label{eq:obj1}
J_{i}(u)=\gamma \int_\Xi |x(T,\xi)|^2\,d\xi - \int_0^T\!\!\int_\Xi
|u(t,\xi)|^{2}\,d\xi\,dt
\end{equation}
over all controls $u\in L_2([0,T]\times \Xi,\erre)$, subject to
the dynamics (\ref{eq:pde}). Here the weight coefficient $\gamma>0$
represents the relative impact of the profit contribution of the final
goodwill level vs. that of the advertising effort.
\end{itemize}
While a quadratic cost of advertising effort is a common and plausible
assumption, the choice of this specific form for the profit
contribution of the final pre-launch goodwill level in (\ref{eq:obj1})
is not an obvious one. In particular, it is not concave, and as such
it can be interpreted as the utility function of a \emph{risk-loving}
firm.  However, in some situations this assumption is reasonable (see
e.g. \cite{grosset} for a related discussion) and can also be justified by
considerations of (local) second-order approximations and of
analytical tractability.

\smallskip

In addition to (P1) we also study a ``targeting'' problem denoted as (P2):
\begin{itemize}
\item[(P2)] minimize the functional 
\begin{equation}
\label{eq:obj2}
J_{h}(u)=\gamma \int_\Xi|x(T,\xi )-k(\xi)|^2\,d\xi
+\int_{0}^{T}\!\!\int_{\Xi }|u(t,\xi )|^{2}\,d\xi \,dt
\end{equation}%
over all controls $u\in L_2([0,T]\times \Xi,\erre)$, subject to the
dynamics (\ref{eq:pde}), where $\gamma>0$ is a weight coefficient as above, 
and $k:\Xi\to\erre$ is the target distribution of goodwill at time $T$.
\end{itemize}
The targeting problem may arise in settings where a firm would like to
evaluate the resources required for establishing a particular goodwill
profile $k(\xi)$ for the product to be launched. In particular, if a
firm is interested in enforcing a uniform good will profile
$k(\xi)=k_0$, setting $\gamma$ appropriately high and solving the
targeting problem repeatedly for several values of $k_0$ will provide
an assessment of necessary advertising spending as a function of
$k_0$. Thus, in essence, the targeting problem serves as a surrogate
penalty-function formulation of the problem in which the final
goodwill profile is set to equal $k(\xi)$.

\section{Reformulating the spatial advertising problem as
a control problem in infinite dimensions}
\label{sec:reform}
Let $X$ be the Hilbert space of square integrable
functions defined on the domain $\Xi$, i.e. $X=L_2(\Xi)$, equipped
with the natural inner product
$$
\langle f,g\rangle :=\int_{\Xi }f(\xi )g(\xi )\,d\xi
$$
and norm 
$$
|f|:=\left( \int_{\Xi }f^{2}(\xi )\,d\xi \right) ^{1/2}.
$$
Denote by $A$ the following linear operator in $X$: 
\begin{equation}
\label{eq:A}
\left\{ 
\begin{array}{l}
Ay=(\Delta-\rho)y, \\[10pt] 
\ds D(A)=\Big\{v\in H^1(\Xi):\, \Delta v\in L_2(\Xi),\,
     \int_\Xi \Delta v\phi\,dx = -\int_\Xi\cp{\nabla v}{\nabla\phi}\,dx\Big\}.
\end{array}\right.  
\end{equation}
Note that, in view of Green's formula, the second condition defining
the domain of $A$ in (\ref{eq:A}) can be seen as a weak formulation of
the Neumann boundary condition $\partial v/\partial n = 0$ on
$\partial\Xi$.

Setting (with a slight abuse of notation) $x(t)=x(t,\cdot)$ and
$u(t)=u(t,\cdot)$, we can equivalently write (\ref{eq:pde}) as an
abstract linear system on the Hilbert space $X$:
\begin{equation}
\label{eq:ade}
\left\{
\begin{array}{l}
\displaystyle \frac{dx}{dt}=Ax(t)+Bu(t), \\[10pt] 
x(0)=x_{0}\in X%
\end{array}%
\right.
\end{equation}
for $t\in \lbrack 0,T]$, with $A:D(A)\subset X\to X$ as in (\ref%
{eq:A}), $u \in L_2([0,T],X)$, 
%$u\in\mathfrak{U}:=\{u:[0,T]\to X,\, u(t,\xi)\in[0,R]
%\;\;(t,\xi)\mathrm{-a.e.}\}$, 
and $B:X\to X$ is the linear
bounded operator defined by
\begin{equation}
\label{eq:B}
B:y(\xi )\mapsto b(\xi )y(\xi ).
\end{equation}
We shall need the following important features of the Neumann
Laplacian, which we collect in the form of a lemma. The proof can be
found e.g. in \cite{RS4} and \cite{vrabie}.
\begin{lemma}\label{propA}
  The linear operator $A$ defined in (\ref{eq:A}) is self-adjoint,
  negative, and generates a strongly continuous positivity-preserving
  semigroup of contractions. Moreover, the Cauchy problem
  (\ref{eq:ade}) admits a unique mild solution $x$ given by the
  variation of constants formula:
\begin{equation}
\label{eq:mild}
x(t)=e^{tA}x_{0}+\int_{0}^{t}e^{(t-s)A}Bu(s)\,ds,
\end{equation}%
where $e^{tA}$ denotes the strongly continuous semigroup generated by $A$.
\end{lemma}

Problem (P) can now be written as 
\begin{equation}
\label{eq:P}
\inf_{u\in \mathfrak{U}}\left( \phi (x(T))+\int_{0}^{T}h_{1}(u(t))\,dt\right)
\end{equation}%
subject to the dynamics (\ref{eq:ade}), where $\phi $, $h_{1}:X\to
\erre$ are defined as
\begin{eqnarray}
\label{eq:phi}
\phi (x) &=&-\int_{\Xi }\phi _{0}(x(\xi ))\,d\xi,\\
h_{1}(u) &=&\int_{\Xi }h_{0}(u(\xi ))\,d\xi, \nonumber
\end{eqnarray}%
and 
$$
\mathfrak{U}=\Big\{u:[0,T]\to X\Big|\;u(\cdot )(\xi)\in [0,R],
\,x(\cdot )(\xi )\geq 0\;\xi \textrm{-a.e.}\Big\}.
$$

Similarly, the objective functionals (\ref{eq:obj1}) and (\ref{eq:obj2}) can
be respectively written as 
$$
J_{i}(u)=-\gamma |x(T)|^{2}+\int_{0}^{T}|u(t)|^{2}\,dt
$$
and 
$$
J_{h}(u)=\gamma |x(T)-k|^{2}+\int_{0}^{T}|u(t)|^{2}\,dt.
$$
The aim is to find an optimal control, i.e. a function
$u_*\in\mathfrak{U}_{ad}$ such that
$$
J(u_*)\leq J(u) \quad \forall u \in \mathfrak{U}_{ad}.
$$
Here $J$ is either $J_{c}$, $J_{i}$, or $J_{h}$, and
$\mathfrak{U}_{ad}$ is the class of admissible controls:
$\mathfrak{U}_{ad}=\mathfrak{U}$ for problem (P), and
$\mathfrak{U}_{ad}=L_2([0,T],X)$ for problems (P1) and (P2). The
pair $(x_*,u_*)$, where $x_*$ is the solution of
(\ref{eq:ade}) with $u\equiv u_*$, is often called the
\emph{optimal pair} for the corresponding optimal control problem.

\begin{rmk}
  Problem (P1) is a linear-quadratic (LQ) optimal control problem with
  indefinite costs, while problem (P2) is an LQ problem with positive
  costs similar to those encountered in the study of target tracking.
  While problem (P2) is always well-posed and always admits an optimal
  control, problem (P1) will be in general only locally well-posed,
  and global well-posedness will follow from additional assumptions on
  the parameters of the problem. For more details on LQ problems with
  indefinite costs in infinite dimensions, see \cite{LY}.
\end{rmk}

\begin{rmk}
The objective functional of problem (P1) could be taken, more
generally, as 
$$
J_{i}(u)=-\left\langle P_{0}x(T),x(T)\right\rangle
+\int_{0}^{T}\left| u(t)\right| ^{2}dt,
$$
with $P_{0}:X\to X$. For example, if the goodwill is more valued in a
region $\Xi_1\subset \Xi$, $P_{0}$ could be an operator such that
$P_{0}=\gamma _{1}I$ on $\Xi \setminus \Xi_1$, with $\gamma
<\gamma _{1} $. Even more generally, one could fix a bounded function
$p:\Xi \to \erre_+$ modeling the importance of goodwill at each
point of $\Xi $, and define $P_{0}=pI$, $I$ being the identity
function on $X$. Note that the above also applies to problem (P2).
\end{rmk}

\section{Solution of the general constrained problem}
The basic idea is to embed the state and control constraints into the
structure of the problem. Observe that the non-negativity of the
intial condition and of the control in (\ref{eq:ade}) implies the
non-negativity of the state variable at any point in time. Using this
observation and assigning infinite cost to the controls that do not
satisfy the constraint $0\leq u\leq R$ we enforce the required state
and control constraints.  Next we show that the problem admits a
solution, i.e. that the optimal control exists, and we write a weak
maximum principle that gives a (abstract) characterization of the
sought optimal advertising policy.  Finally we discuss two special
cases of the problem for which we can find explicit solutions from the
abstract maximum principle.

Let us define $h:X\to \erre \cup \{+\infty\}$ as follows: 
\begin{equation}
\label{eq:h}
h(u)=\left\{ 
\begin{array}{ll}
h_{1}(u), & u\in \lbrack 0,R] \textrm{ a.e.}, \\[8pt] 
+\infty , & \textrm{otherwise}.%
\end{array}%
\right.
\end{equation}

\begin{prop}\label{prop:exist}
Assume that
\begin{enumerate}
\renewcommand{\labelenumi}{{\rm (\roman{enumi})}}
\item $h_{0}:\mathbb{R} \to \mathbb{R}$ is convex and 
lower semicontinuous;
\item $\phi_{0}$ is continuous and concave;
\item $x_{0}$ is nonnegative and $x_{0}\in H_{0}^{1}\left( \Xi
\right) $.
\end{enumerate}
Then there exists at least one optimal pair $(x_*,u_*)$ for problem
(P), with $x_*\in C([0,T],X) $ and $u_*\in\mathfrak{U}$.
\end{prop}
The proof of this proposition is essentially classical, but we report
it for the reader's convenience and for non-specialists in infinite
dimensional control theory.
\begin{proof}
First note that $\mathfrak{U}\subset L_1([0,T],X) \cap L_2([0,T],X)$
because $\Xi$ is a bounded set in $\erre^2$.
For every $u\in L_1([0,T],X)$ the state system (\ref{eq:ade}) admits
a unique mild solution $x^u\in C([0,T],X)$. Define the function
\begin{eqnarray*}
\Phi : L_2([0,T],X) &\to& \mathbb{R} \\
u &\mapsto &\int_0^T h(u(t))\,dt + \phi(x^u(T)).
\end{eqnarray*}
Problem (P) can be equivalently written as
$$
\inf_{u\in L_2\left( \left[ 0,T\right] ,X\right) }\Phi \left( u\right) .
$$
Let us now show that $\Phi$ is convex and lower semicontinous: for
$\lambda\in[0,1]$ and $u_1$, $u_2 \in L_2([0,T],X)$, the convexity of
$h$ implies
\begin{equation}
  \label{eq:conv1}
  \int_0^T h(\lambda u_1(t) + (1-\lambda)u_2(t))\,dt \leq
     \lambda \int_0^T u_1(t)\,dt + (1-\lambda)\int_0^T u_2(t)\,dt.
\end{equation}
One also has
\begin{eqnarray*}
x^{\lambda u_1 + (1-\lambda)u_2}(T) &=& e^{-AT}x_0 
  + \int_0^T e^{-A(T-t)}B(\lambda u_1(t) + (1-\lambda)u_2(t))\,dt \\
&=& \lambda e^{-AT}x_0 + (1-\lambda)e^{-AT}x_0 \\
  && \quad  + \lambda \int_0^T e^{-A(T-t)}Bu_1(t)\,dt + 
      (1-\lambda)\int_0^T e^{-A(T-t)}Bu_2(t)\,dt \\
&=& \lambda \Big( e^{-AT}x_0 + \int_0^T e^{-A(T-t)}Bu_1(t)\,dt \Big) \\
  && \quad + (1-\lambda) \Big( e^{-AT}x_0 + \int_0^T e^{-A(T-t)}Bu_2(t)\,dt
\Big) \\
&=& \lambda x^{u_1}(T) + (1-\lambda) x^{u_2}(T),
\end{eqnarray*}
hence, by the convexity of $\phi$,
\begin{equation}
  \label{eq:conv2}
  \phi\Big(x^{\lambda u_1 + (1-\lambda)u_2}(T)\Big) \leq
     \lambda\phi(x^{u_1}(T)) + (1-\lambda) \phi(x^{u_2}(T)).
\end{equation}
Inequalities (\ref{eq:conv1}) and (\ref{eq:conv2}) imply
$$
\Phi(\lambda u_1 + (1-\lambda)u_2) \leq \lambda\Phi(u_1) + (1-\lambda)\Phi(u_2),
$$
i.e. $\Phi$ is convex. In order to prove lower semicontinuity, let
$(u_n)_{n\geq 0}$ be a sequence which converges strongly to $u$ in
$L_2([0,T]\times\Xi)$. Then $u_n \to u$ in measure, and in particular
there exists a subsequence $u_{n_k}$ such that $u_{n_k} \to u$
$(t,\xi)$-a.e.. By well-known properties of convex functions, $h$ is
bounded from below by an affine function, i.e. there exists $a$,
$b\in\erre$ such that $h(x)\geq ax+b$, hence
$$ 
\lim_{n_k\to\infty} h(u_{n_k})-au_{n_k}-b = h(u)-au-b
\quad (t,\xi)\textrm{-a.e.},
$$
and by Fatou's lemma
$$
\liminf_{n_k\to\infty} \int_{[0,T]\times X} (h(u_{n_k})-au_{n_k}-b)\,d\xi\,dt
\geq \int_{[0,T]\times X} (h(u)-au-b) \,d\xi\,dt,
$$
i.e.
$$
\liminf_{n_k\to\infty} \int_{[0,T]\times X} (h(u_{n_k}-h(u))\,d\xi\,dt
+ a \liminf_{n_k\to\infty} \int_{[0,T]\times X} (u-u_{n_k})\,d\xi\,dt \geq 0.
$$
But
\begin{eqnarray*}
\left| 
\liminf_{n_k\to\infty} \int_{[0,T]\times X} (u-u_{n_k})\,d\xi\,dt
\right|
&\leq&
\liminf_{n_k\to\infty} \int_{[0,T]\times X} |u-u_{n_k}|\,d\xi\,dt \\
&\leq& (T\,\mathrm{vol}(\Xi))^{1/2} \liminf_{n_k\to\infty} 
\left[\int_{[0,T]\times X} |u-u_{n_k}|^2\,d\xi\,dt\right]^{1/2} \\
&=& 0,
\end{eqnarray*}
hence
$$
\liminf_{n_k\to\infty} \int_0^T h(u_{n_k}(t))\,dt \geq 
\int_0^T h(u(t))\,dt.
$$
Completely similar arguments show that $\liminf_{n_k\to\infty}
\phi(x^{u_{n_k}}(T) \geq \phi(x^u(T))$, thus we have proved that
$\Phi$ is lower semicontinuous.\\
Moreover, one has $\lim_{|u|\to\infty} \Phi(u)=+\infty$.
Let now $a\in\erre$ be any (fixed) number. Set $E=\{x\in X:\Phi(x)\leq
a\}$. Then one has
\begin{equation}
  \label{eq:E}
\inf_{x\in X}\Phi(x) = \inf_{x\in E}\Phi(x).  
\end{equation}
The lower semicontinuity of $\Phi$ implies that the level set $E$ is
closed, hence weakly closed. Moreover, since $\lim_{|u|\to\infty}
\Phi(u)=+\infty$, one also has that $E$ is bounded. Therefore $E$ is
also compact in the weak topology of $X$. By a well known result,
every lower semicontinuous function on a compact subset of a
topological space attains its infimum. The existence of a minimizer
now follows immediately from (\ref{eq:E}).\\
Let $u_*$ be the minimizer of $\Phi $. Then it is clear that it
must be $u_*\in\mathfrak{U}$. By lemma \ref{propA} we also have that
$e^{tA}$ is a positivity preserving semigroup.  If $u_*\geq 0$ a.e.,
then one also has $Bu_*\geq 0$ a.e. by the assumptions on $\xi\mapsto
b(\xi)$. This implies, together with (\ref{eq:mild}), assumption
(iii), and the positivity preserving property of $e^{tA}$, that
$x_*\geq 0$ a.e., which concludes the proof.
\end{proof}
Note that this proof of existence of an optimal pair can be adapted to
allow for more general type of constraints. Consider for instance the
problem 
$$
\sup_{u}\phi (x^{u}(T)),
$$
where $u$ satisfies the budget constraint 
$$
\int_{0}^{T}\int_{\Xi }u^{2}(t,\xi )\,d\xi \leq M,
$$
with $M$ a fixed positive number, and $\phi $ is a concave function. Here we
have followed \cite{muller} in using the quadratic form for the cost as a
function of advertising effort. Then we define the set of admissible
controls ${\cal U}\subset L_2([0,T],X)$ as 
$$
{\cal U}=\left\{ u: \; \int_{0}^{T}\int_{\Xi }u^{2}(t,\xi )\,d\xi
\leq M\right\} ,
$$
and the functional $\Phi :L_2([0,T],X)\to \erre$ as 
$$
\Phi :u\mapsto -\phi (x^{u}(T)).
$$
Now we can reformulate the problem as 
$$
\inf_{u\in\mathcal{U}} \Phi (u).
$$
In complete analogy to the proof of proposition \ref{prop:exist}, one
could show that the set $\mathcal{U}$ is convex, closed, and bounded
in $L_2([0,T],X)$, that the function $\Phi$ is convex lower
semicontinuous in $L_2([0,T],X)$, and that these conditions are enough
to ensure that an optimal pair $(x_*,u_*)$ exists, with
$u_*\in\mathcal{U}$. %as follows by an application of Proposition 1.2
%in Chapter 2 of \cite{ET}.

Optimal pairs in many cases can be characterized by the first-order
necessary optimality conditions, through the introduction of a Lagrange
multiplier. In particular, one has to solve 
$$
\inf_{u\in L_2([0,T],X)}\left[ -\phi (x(T))+\lambda
\int_{0}^{T}u^{2}(t)\,dt\right] .
$$
Let $u_{\lambda }$ be the optimal control for this problem. Then $u_*$
is given by a $u_{\lambda }$ such that 
$$
\int_{0}^{T}u_{\lambda }^{2}(t)\,dt=M.
$$

Once the existence of an optimal solution is established, we can characterize
optimal pairs using the maximum principle in Hilbert spaces.

\begin{prop}
  Let $(x_*,u_*)$ be an optimal pair for
  problem (P). Then there exists $p\in
  C\left([0,T],X\right)$ such that $x_*,$ $u_*,$ $p$
  satisfy the following two-point boundary value problem:
\begin{equation}
\label{eq.maxpr}
\left\{ 
\begin{array}{l}
\displaystyle x'_*=Ax_*+Bu_*, \\[4pt] 
\displaystyle p'+Ap=0, \\[4pt]
x_*(0)=x_0, \\[4pt] 
p(T) \in -\partial \phi( x_*(T)), \\[4pt] 
u_*(t) \in \partial h^*(Bp(t)),
\end{array}
\right.
\end{equation}
where $h^*$ is the conjugate of $h$, and $\cdot'$ stands for
differentiation with respect to time.
\end{prop}
\begin{proof}
  In order to apply the sufficient and necessary optimality system of
  Theorem 4.2.1 of \cite{barbu95}, we need to check some conditions.
  In particular, $A$ is the generator of a strongly continuous
  semigroup (lemma \ref{propA}), and $B$ is a linear bounded operator,
  as follows from its definition (\ref{eq:B}) and the boundedness of
  $b$. The continuity and convexity of $\phi$ follow immediately by
  its definition (\ref{eq:phi}), and similarly for the convexity and
  lower semicontinuity of $h$, cf.  (\ref{eq:h}). Therefore the
  results of the theorem hold, from which we get our optimality
  conditions (\ref{eq.maxpr}) by obvious modifications of those in
  \cite{barbu95}.
\end{proof}
In the next section we consider two special cases of the general
spatial advertising problem for which we can provide more specific
characterizations of the optimal advertising policies.

\subsection{Two cases with explicit solutions}
Deriving explicit expressions for the optimal pair from the maximum
principle is in general a very challenging task. On the other hand,
for some specific choices of the objective function, these closed-form
solutions can be obtained. Consider, for instance, the case where $\phi
(x)=-x $, and $h(u)=|u|^{2}/2$, for $u\in \lbrack 0,R]$ a.e., and
$h(u)=+\infty $ otherwise. Then one has $-\partial \phi =1$, and
$$
h^{\ast }(\zeta )=\left\{ 
\begin{array}{ll}
0, & \zeta <0, \\[4pt] 
\zeta ^{2}/2, & \zeta \in \lbrack 0,R], \\[4pt] 
\zeta R-R^{2}, & \zeta >R.%
\end{array}%
\right.
$$
Therefore the maximum principle (\ref{eq.maxpr}) yields
$$
p'+Ap=0,\qquad p(T)=1,
$$
hence, by well known properties of the heat equation with Neumann
boundary conditions,
$$
p(t,\xi) = e^{-\rho(T-t)}, \quad (t,\xi) \in [0,T]\times \Xi.
$$
From this representation of the dual arc $p$ we immediately obtain
that the optimal strategy is given by 
\begin{equation}
\label{eq:capri}
u_*(t,\xi)=b(\xi)e^{-\rho(T-t)}\wedge R,
\qquad
\textrm{a.e.-}(t,\xi) \in [0,T]\times \Xi.
\end{equation}
In particular, $u_*$ is always nonnegative and is zero where the
advertising effectiveness $b$ is zero.  Finally, the optimal
trajectory $x_*$ is given by
$$
x_*(t)=e^{tA}x_{0}+\int_{0}^{t}e^{(t-s)A}Bu_*(s)\,ds. 
$$
An equivalent expression for the optimal trajectory that is easier to
implement numerically can be given by projecting on a basis of
$L_2(\Xi)$ of eigenvectors of $A$. In fact, it is well known (see e.g.
\cite{agmon}) that there exists a complete orthonormal (ONC) system
$(e_k)_{k\in\mathbb{N}}$ in $L_2(\Xi)$ and a sequence of positive
numbers $(\lambda_k)_{k\in\mathbb{N}} \uparrow +\infty$ such that
$$
Ae_k = -(\lambda_k+\rho)e_{k},\quad \forall k\in \N.
$$
Then we have
$$
x_*(t,\xi) = \sum_{k\in\N} x_*^k(t)e_k(\xi),
$$
a.e. in $\Xi_T$, where
$$
x_*^k(t) = e^{-(\lambda_k+\rho)t}x^k(0)
+ \int_0^t e^{-(\lambda_k+\rho)(t-s)}\tilde{u}^k(s,\xi)\,ds,
$$
and $x^k(0)=\cp{x_0}{e_k}_{L_2(\Xi)}$,
$\tilde{u}^k(t)=\cp{b(\cdot)u(t,\cdot)}{e_k}_{L_2(\Xi)}$.

An analogous procedure allows one to obtain the optimal advertising policy
and the optimal state with linear cost, i.e. with $h(u)=u$ for $u\in \lbrack
0,R]$ a.e., $h(u)=+\infty $ otherwise. In particular 
\begin{equation}
\partial h^{\ast }(p)=\left\{ 
\begin{array}{ll}
0, & p<1, \\ 
R, & p>1,%
\end{array}%
\right.  \label{eq:hstar}
\end{equation}%
hence 
\begin{equation}
\label{eq:cciosa}
u_*(t,\xi )=\left\{ 
\begin{array}{ll}
0, & b(\xi )p(t,\xi )<1, \\ 
R, & b(\xi )p(t,\xi )> 1,%
\end{array}%
\right.
\end{equation}%
which solves the problem.

Unfortunately we were unable to obtain closed form expressions for the
modified version of problem (P2) with linearized cost of control.
While the general existence criterion can still be applied, the
two-point boundary value problem that characterizes the optimal pair
cannot be solved explicitly. In fact, the maximum principle yields
\begin{equation}
\label{eq:quattro}
\left\{
\begin{array}{l}
\displaystyle x'_*=Ax_* + Bu_*, \\[4pt] 
\displaystyle p'+Ap=0, \\[4pt] 
x_*(0)=x_0, \\[4pt] 
p(T)=-2(x_*(T)-k), \\[4pt] 
u_*(t) \in \partial h^{\ast }(Bp(t)),
\end{array}
\right.
\end{equation}
where $\partial h^{\ast }$ is given by (\ref{eq:hstar}).
\begin{rmk}
  Note that in the absence of the
  diffusion term, i.e. for $A=-\rho$, the optimal strategies
  (\ref{eq:capri}) and (\ref{eq:cciosa}) would remain unchanged
  (because the dual arc $p$ would not change), even though the
  corresponding optimal trajectory do not coincide, as a simple
  calculation reveals. However, this if of course not a general
  phenomenon, but just a consequence of the equality $\partial\phi=-1$
  in the above examples. For instance, such simplification does not
  happen in (\ref{eq:quattro}).
\end{rmk}

In the next section we employ the dynamic programming principle to
give a characterization of optimal strategies in feedback form. Even
though we cannot obtain, in general, completely explicit
representation for the optimal strategy, we shall at least indicate
some schemes to approximate the value function associated to the
problem.

\section{Optimal feedback strategies}
Let us define the value function
\begin{equation}
\label{eq:value}
V(t,y) = \inf_{u\in\mathfrak{U}}
\Big( \phi(x(T)) + \int_t^T h(u(s))\,ds \Big),
\end{equation}
subject to
\begin{equation}
\label{eq:ddr}
x' = Ax + Bu,
\quad s\in [t,T],
\quad x(t)=y.
\end{equation}
A mild solution of (\ref{eq:ddr}) will be denoted by $x^{y,u}$.
Note that, strictly speaking, the value function of problem (P) is
$-V(t,y)$. However, for consistency with the previous sections, we
prefer to study the minimization problem (\ref{eq:value}), which is
clearly equivalent to problem (P).

The existence of an optimal pair $(x_*,u_*)$ for (\ref{eq:value}) is
guaranteed by proposition \ref{prop:exist}. We shall now see that the
dynamic programming approach allows one to obtain a characterization
of $(x_*,u_*)$ in terms of the solution of an associated
Hamilton-Jacobi equation.

We first prove two simple but important qualitative properties of the
value function.
\begin{prop}
The value function $V(t,y)$ is convex and decreasing in $y$.
\end{prop}
\begin{proof}
Let $y_1$, $y_2\in X$, and assume that $u_1$, $u_2$ are optimal for
$V(t,y_1)$ and $V(t,y_2)$, respectively.
Then we have, for any $\lambda\in[0,1]$,
\begin{eqnarray*}
\lambda V(t,y_1) + (1-\lambda)V(t,y_2) &=&
\lambda \phi(x^{y_1,u_1}(T))
+ (1-\lambda) \phi(x^{y_2,u_2}(T))\\
&& + \int_t^T \lambda h(u_1(s))\,ds
   + \int_t^T (1-\lambda) h(u_2(s))\,ds\\
&\geq& \phi\big(\lambda x^{y_1,u_1}(T) + (1-\lambda)x^{y_2,u_2}(T)\big)\\
&& + 
\int_t^T [\lambda h(u_1(s)) + (1-\lambda) h(u_2(s))]\,ds,\\
\end{eqnarray*}
where the inequality holds because $\phi$ is convex.
A simple calculation reveals that
$$
\lambda x^{y_1,u_1}(T) + (1-\lambda)x^{y_2,u_2}(T)
= x^{\lambda y_1 + (1-\lambda)y_2,u}(T),
$$
where $u=\lambda u_1 + (1-\lambda)u_2 \in \mathfrak{U}$. By the
convexity of the cost function $h$ we then have
\begin{eqnarray*}
\lambda V(t,y_1) + (1-\lambda)V(t,y_2) &\geq&
\phi(x^{\lambda y_1 + (1-\lambda)y_2,u}(T)) + \int_t^T h(u(s))\\
&\geq& \inf_{u\in\mathfrak{U}} \Big(
       \phi(x^{\lambda y_1 + (1-\lambda)y_2,u}(T)) 
       + \int_t^T h(u(s)) \Big)\\
&=& V(t,\lambda y_1 + (1-\lambda)y_2).
\end{eqnarray*}
Let us now show that $V(t,y)$ is increasing in $y$: let $u_*$ be an
optimal control, so that
$$
V(t,y) = \phi(x^{y,u_*}(T)) + \int_t^T h(u_*(s))\,ds,
$$
and take $\bar y \geq y$. Then we have
\begin{eqnarray*}
J(t,\bar y,u_*) &=& \phi(x^{\bar y,u_*}(T)) + \int_t^T h(u_*(s))\,ds\\
&=& \phi\Big(e^{(T-t)A}\bar y + \int_t^T e^{(T-t-s)A}Bu_*(s)\,ds\Big)
    + \int_t^T h(u_*(s))\,ds\\
&\leq& \phi\Big(e^{(T-t)A}y + \int_t^T e^{(T-t-s)A}Bu_*(s)\,ds\Big)
    + \int_t^T h(u_*(s))\,ds\\
&=& V(t,y),
\end{eqnarray*}
where the inequality holds because the semigroup generated by $A$ is
positivity preserving (cf. lemma \ref{propA}), $b(\xi)\geq 0$ $\xi$-a.e.,
and $\phi_0$ is increasing. The proof is completed observing that
$V(t,\bar y) \geq J(t,\bar y,u_*)$.
\end{proof}

A combination of the dynamic programming principle and of the abstract
maximum principle allows one to write a feedback representation of the
optimal strategy in terms of the value function $V$.
\begin{prop}
Let $(x_*,u_*)$ an optimal pair for (\ref{eq:value}), and $V:[0,T]\times X \to
\erre$ the corresponding value function. Then the following inclusion holds:
$$
u_*(s) \in \partial h^*(-B^* \partial V(s,x_*(s))),
$$
where $\partial V$ stands for the subdifferential of the function
$y \mapsto V(s,y)$.
\end{prop}
\begin{proof}
Let $v \in L_1([0,T],X)$, and let $z$ be a (mild) solution of
$z'=Az+Bv$. Then we can write
\begin{eqnarray*}
  \lefteqn{\cp{x_*(T)-z(T)}{p(T)} - \cp{x_*(s)-z(s)}{p(s)}}\\
&=& \int_s^T \cp{x'_*(r)-z'(r)}{p(r)}\,dr 
    + \int_s^T \cp{x_*(r)-z(r)}{p'(r)}\,dr\\
&=& \int_s^T \cp{A(x(r)-z(r))}{p(r)}\,dr 
    + \int_s^T \cp{B(u_*(r)-v(r))}{p(r)}\,dr\\
& & + \int_s^T \cp{x_*(r)-z(r)}{p'(r)}\,dr\\
&=& \int_s^T \cp{B(u_*(r)-v(r))}{p(r)}\,dr,
\end{eqnarray*}
where the second equality follows by definitions of $y$ and $z$, and
the third equality follows by the identity $p'\in -A^*p$. By
definition of subdifferential we have
$$
\cp{B(u_*(r)-v(r))}{p(r)} \geq
h(u_*(r))-h(v(r)),
$$
because $u_*(t)\in \partial h^*(B^*p(t))$ implies
$B^*p(t)\in\partial h(u_*(t))$.
Recalling that $p(T) \in -\partial\phi(x_*(T))$
we obtain
$$
- \cp{x_*(s)-z(s)}{p(s)} \geq
\phi(x_*(T)) + \int_s^T h(u_*(r))\,dr
-\varphi(z(T)) - \int_s^T h(v(r))\,dr.
$$
The dynamic programming principle now implies
$$
V(s,x_*(s)) = 
\phi(x_*(T)) + \int_s^T h(u_*(r))\,dr,
$$
which together with the previous inequality yields $-p(s) \in \partial
V(s,x_*(s))$. Recalling that $u_*(s) \in \partial h^*(B^*p(s))$ we
finally obtain
$$
u_*(s) \in \partial h^*(-B^*\partial V(s,x_*(s))).
\qedhere
$$
\end{proof}

One would expect that the value function $V$ solves, in a suitable
sense, the Hamilton-Jacobi equation associated to problem
(\ref{eq:value}), which can be written as
\begin{equation}
\label{eq:apri}
v_s(s,y) + \cp{Ay}{v_y(s,y)} + F(v_y(s,y))=0,
\quad v(T,y)=\phi(y), 
\end{equation}
for $s\in[t,T]$ and $y\in D(A)$, where $F(q):=-h^*(-B^*q)$.  In fact
(see \cite{BDP-ann}), assuming that $F \in C^1(H)$, $V$ satisfies
(\ref{eq:apri}) in the sense that there exists $\eta(s,y)\in\partial
V(s,y)$ for all $(s,y)\in [0,T]\times X$ such that
$$
\left\{
\begin{array}{l}
\ds v_s(s,y) + \cp{Ay}{\eta(s,y)} + F(\eta(s,y))=0,
\quad \textrm{a.e.\ }t\in[0,T],\;\; y\in D(A),\\[4pt]
v(T,y)=\phi(y).
\end{array}\right.
$$
In order to approximate the solution of this terminal value problem,
one could apply the Crandall-Liggett theorem, which ensures that $V$
is the limit of a sequence of implicit difference approximations
(time-discretization), under some extra assumptions. For instance, if
$B=I$ and $h(z)=|z|^2$, then $V_\varepsilon(t,\cdot) \to V(t,\cdot)$
uniformly in $t$ as $\varepsilon \to 0$, where
$V_\varepsilon(t,\cdot)$ solves
$$
\left\{
\begin{array}{ll}
\ds \frac{1}{\varepsilon} (V_\varepsilon(s,y)-V_\varepsilon(s-\varepsilon,y))
+ \cp{Ay}{\nabla V_\varepsilon(s,y)}
+ F(\nabla V_\varepsilon(s,y)) = 0,
& s\in]0,T],\;\; y\in D(A),\\[8pt]
V_\varepsilon(s,y)=\phi(y),
&s=0,\;\; y\in D(A).
\end{array}\right.
$$
A formal justification can be given following \cite{BDP-ann}, sec. 4.

\section{Quadratic cost and reward functions}
In this section we solve the optimal control problem (P1) through the
dynamic programming approach, that is, first we solve the associated
operator Riccati equations, and then we show that the optimal control
$u_*$ can be written as a linear feedback of the optimal trajectory
$x_*$.
Let us recall that problem (P1) amounts to minimizing the objective
functional
$$
J_i(u) = -\gamma |x(T)|^2 + \int_0^T |u(t)|^2\,dt,
$$
subject to the dynamics $x'=Ax+Bu$, $x(0)=x_0$.

In this and the following section we also assume $B=I$, for
simplicity. The analysis of the more general case $B\in \mathcal{L}(X)$
follows essentially the same steps.

An abstract solution of problem (P1) is given by the following proposition.
\begin{prop}
Suppose that the parameters of problem (P1) satisfy the
inequality
$$
1-\frac{\gamma }{2\rho }\left( 1-e^{-2\rho T}\right) >0.
$$
Then for any $x_{0}\in X$, $J_{i}$ admits a unique minimizer
$u_*$ given by
$$
u_*=\gamma (I-\gamma L_T^*L_T)^{-1}L_T^*e^{-AT} x_0,
$$
with $L_T$ defined in (\ref{eq:LT}).
%Moreover, the value
%function is a bounded bilinear form on $X$ explicitly given by
%\begin{equation}
%V\left( y\right) =\inf_{u}J\left( y;u\right) =\left\langle \Gamma
%y,y\right\rangle -\left\langle \Psi ^{-1}\Theta y,\Theta y\right\rangle
%\end{equation}%
%for all $y\in X$, where $\Gamma =e^{AT}P_{0}e^{AT}=-\gamma e^{2AT}$.
\end{prop}
\begin{proof}
Let us define the following linear operator:
\begin{eqnarray}
L_T: L_2([0,T],X) &\to& X \nonumber \\
u &\mapsto &\int_{0}^{T}e^{(T-s)A}u(s)\,ds.
\label{eq:LT}
\end{eqnarray}
Then we have $x(T) = e^{-AT}x_0 + L_Tu$, and, denoting by $\|\cdot\|$
the norm in $L_2([0,T],X)$,
\begin{eqnarray}
J_i(u) &=& \cp{P_0x(T)}{x(T)} + \|u\|^2 \nonumber\\
&=& \cp{P_0(e^{-AT}x_0+L_Tu)}{e^{-AT}x_0+L_Tu} + \|u\|^2 \nonumber\\
&=& \cp{e^{-A^*T}P_0e^{-AT}x_0}{x_0} + 2\cp{P_0e^{-AT}x_0}{L_Tu} +
\cp{P_0L_Tu}{L_Tu} + \|u\|^2 \nonumber\\
&=& \cp{e^{-AT}P_0e^{-AT}x_0}{x_0} + 2\cp{L_T^*P_0e^{-AT}x_0}{u} +
\cp{(I+L_T^*P_0L_T)u}{u}.\label{eq:J-conv}
\end{eqnarray}
By (\ref{eq:J-conv}) it immediately follows that if $I+L_T^*P_0L_T
>\!\!\!>0$, then the functional $J_i:L_2([0,T],X) \to \erre$ is convex,
lower semi-continuous, and such that
$\lim_{|u|\to\infty} J_i(u) = +\infty$, hence it attains its minimum in
$L_2([0,T],X)$. Moreover, again as a consequence of (\ref{eq:J-conv}),
if $u_*$ is a minimizer of $J$, then one must have
\begin{equation}
\label{eq:LQi-ol}
(I+L_T^*P_0L_T)u_* + L_T^*P_0e^{-AT}x_0 = 0.
\end{equation}
If $\Psi:=(I+L_T^*P_0L_T)$ is an homeomorphism of $L_2([0,T],X)$ (i.e.
if $\Psi^{-1}$ is also a linear continuous operator on
$L_2([0,T],X)$), then we deduce from (\ref{eq:LQi-ol}) that the unique
optimal control $u_*$ is given by
$$
u_* = -\Psi^{-1}L_T^*P_0e^{-AT}x_0.
$$
From the above discussion, in order to prove the proposition, we
need to show that under the stated hypotheses $\Psi$ is uniformly
positive definite, and $\Psi^{-1}$ is continuous. In order to prove
the former, write
\begin{eqnarray*}
\cp{\Psi v}{v} &=& \|v\|^{2}+\cp{L_T^*P_0L_Tv}{v} =
\|v\|^{2} -\gamma \cp{L_T^*L_Tv}{v} \\
&=& \|v\|^2 -\gamma \cp{L_Tv}{L_Tv} \\
&=& \|v\|^{2}-\gamma |L_Tv|^2,
\end{eqnarray*}
where we have used the fact that $P_0=-\gamma I$. We also have
\begin{eqnarray*}
\left| L_{T}v \right| &=& \left| \int_{0}^{T}e^{-\rho(T-s)}
e^{\Delta(T-s)}v(s)\,ds \right| \\
&\leq &\int_{0}^{T}e^{-\rho \left( T-s\right) }\left\vert v\left( s
\right) \right\vert ds \\
&\leq & \left[ \int_0^T e^{-2\rho(T-s)}\,ds\right]^{1/2}
\left[ \int_{0}^{T} |v(s)|^2\,ds\right]^{1/2} \\
&=&\left( C_{\rho ,T}^{{}}\right) ^{1/2} \| v \|.
\end{eqnarray*}
We have used (in order) a standard estimate, the contractivity of the heat
semigroup, the Cauchy-Schwarz inequality, and Fubini's theorem for positive
integrands. By 
$$
C_{\rho ,T}=\int\limits_{0}^{T}e^{-2\rho \left( T-s\right) }ds=\frac{1}{%
2\rho }\left( 1-e^{-2\rho T}\right)
$$
one has $\cp{\Psi v}{v} \geq \varepsilon \|v\|^2$, with $\varepsilon=
1-\frac{\gamma }{2\rho }\left( 1-e^{-2\rho T}\right)$.\\
In order to prove the continuity of $\Psi$, it is enough to recall the
following simple fact: if $\cp{\Psi v}{v}\geq \varepsilon\|v\|^2$ for
all $v\in D(\Psi)$ and $\Psi$ is self-adjoint, then $\|\Psi^{-1}\|
\leq \varepsilon^{-1}$. But $\Psi=I-\gamma L_T^*L_T$ is clearly
self-adjoint, and the proof is finished.
\end{proof}
This optimal control suffers of two major drawbacks: it is of the open-loop
type, and it is difficult to write explicitly (in particular, the
computation of $\Psi^{-1}$ does not seem to be a straightforward task).
Appealing to the dynamic programming principle, one can obtain a more
explicit feedback characterization of the optimal policy:
\begin{prop}
If
\begin{equation}
1-\frac{\gamma }{2\rho }\left( 1-e^{-2\rho T}\right) >0,
\label{eq:ass-indef}
\end{equation}%
then there exists $P\in C\left( \left[ 0,T\right] ;{\cal L}\left(
X\right) \right) $, $P(t)=P(t)^{\ast }$ {\it for all} $t\in \left[ 0,T\right]
$ which solves the operator Riccati equation
\begin{equation}
\left\{ 
\begin{array}{l}
\displaystyle P'=2AP-PBB^*P, \\ 
P(0)=P_{0}=-\gamma I.%
\end{array}%
\right.  \label{eq:opr1}
\end{equation}%
Moreover, the optimal control is given by $u_*(t)=-B^*P\left(
T-t\right) x_*(t)$, where $x_*(t)$ is the unique (mild)
solution of the closed loop equation 
\begin{equation}
\label{eq:cle1}
\left\{ 
\begin{array}{l}
\displaystyle x'=Ax(t)-BB^*P(T-t)x(t), \\ 
x(0)=x_{0},%
\end{array}%
\right.
\end{equation}%
and the value function can be written as $V(t,y)=J_{i}(u_*)=\cp{P(T-t)y}{y}$.
\end{prop}
\begin{proof}
  The proof follows the classical scheme based on solving the Riccati
  equation and applying a verification theorem (see e.g. \cite{rep2}
  or \cite{DP-ICTP}), but one needs to take into account that $P_0$ is
  \emph{negative} in our case. However, the coercivity of $J_i$ on
  $L_2([0,T],X)$ implied by (\ref{eq:ass-indef}), ``saves'' the
  argument (see also theorem 9.4.3 of \cite{LY}), and in fact already
  gives existence and uniqueness of the optimal pair $(u_*,x_*)$. Here
  we limit ourselves to sketch the proof only, pointing out where the
  essential difference is with respect to the usual case. The first
  step consists in proving that the Riccati equation (\ref{eq:opr1})
  admits a local solution in the space $C([0,\tau];\mathcal{L}(X))$,
  $\tau<T$ (in the usual case one has $\Sigma(X)$, the space of
  continuous linear symmetric operators, instead of $\mathcal{L}(X)$),
  that global uniqueness holds, and finally that local solutions can
  be extended to the whole interval $[0,T]$. The second step is to
  show that the closed-loop equation (\ref{eq:cle1}) admits a unique
  mild solution $x\in C([0,T],X)$, which follows as in the classical
  case.  Again following the classical case, in the third step one
  obtains the representation
$$
J_i(u) = \int_0^T |u(t)+B^*P(T-t)x(t)|^2\,dt + \cp{P(T)x_0}{x_0},
$$
hence $J_i(u)\geq\cp{P(T)x}{x}$, for any $u\in L_2([0,T],X)$.
Therefore, setting $u_*=-B^*P(T-t)x_*(t)$, with $x_*$ the solution of
the closed-loop equation (\ref{eq:cle1}), the proposition is proved.
\end{proof}

\subsection{An example with explicit solution}
In this section we assume that $B=I$, for simplicity, and we show that
the Riccati equation of Proposition 4 can actually be solved
explicitly in this case. In fact, it is well known (see e.g.
\cite{agmon}) that there exists a complete orthonormal (ONC) system
$(e_k)_{k\in\mathbb{N}}$ in $L_2(\Xi)$ and a sequence of positive
numbers $(\lambda_k)_{k\in\mathbb{N}} \uparrow +\infty$ such that
$$
Ae_{k}=-\lambda _{k}e_{k},\quad \forall k\in \N.
$$
Then we can ``project'' the Riccati equation on this ONC system,
obtaining the infinite set of Cauchy problems
\begin{equation}
\label{eq:pk}
\left\{ 
\begin{array}{l}
\displaystyle \frac{dp_{k}}{dt}=-2\lambda _{k}p_{k}(t)-p_{k}^{2}(t), \\[8pt] 
p_{k}(0)=-\gamma ,%
\end{array}%
\right.  
\end{equation}%
where $p_{k}(\cdot ):=P(\cdot )e_{k}$. It is immediate that $q_{k}(t)\equiv
-2\lambda _{k}$ is a particular solution for the $k$-th problem. Then set 
$$
z_k(t)=\frac{1}{p_{k}(t)-q(t)}=\frac{1}{p_k(t)+2\lambda_k}.
$$
One easily obtains that $z_{k}$ satisfies the linear equation 
$$
\frac{dz_{k}}{dt}=-2\lambda _{k}z_{k}(t)+1,
$$
whose general solution is given by 
$$
z_{k}(t)={\frac{1}{2\lambda _{k}}}+C_{k}e^{-2\lambda _{k}t}.
$$
This in turn implies that the general solution for (\ref{eq:pk}) is given by 
\begin{equation}
\label{eq:pks}
p_{k}(t)=-2\lambda _{k}+\frac{1}{(2\lambda _{k})^{-1}+C_{k}e^{-2\lambda
_{k}t}},  
\end{equation}
and by the initial condition
$C_{k}=\frac{\gamma}{2\lambda_k(2\lambda_k-\gamma)}$.

An explicit expression for the optimal trajectory can also be
obtained, again using a projection on the orthonormal system
$(e_k)_{k\in\mathbb{N}}$. In particular, one has
$$
\frac{dx^k}{dt}=-\lambda _{k} x^k(t)+p_{k}(T-t)x^k(t)
$$
for each $k$, hence 
\begin{equation}
\label{eq:traj2-sol}
x^k_*(t)=x^k(0)e^{-\lambda _{k}t}e^{\int_{0}^{t}p_{k}(T-s)\,ds},
\end{equation}
and 
$$
x_*(t,\xi)=\sum_{k=0}^{\infty }x^k_*(t)e_{k}(\xi ).
$$

We can now write explicitly, in terms of the basis $(e_{k})_{k\in \N}$,
the optimal distributed strategy. Namely, 
$$
u^k_*(t)=p_{k}(T-t)x^k_*(t),
$$
with $p_k$ as in (\ref{eq:pk}), and $x^k_*$ is given by
(\ref{eq:traj2-sol}), hence 
$$
u_*(t,\xi)=\sum_{k=1}^{\infty }p_{k}(T-t)x^k_*(t)e_{k}(\xi).
$$

In general it is not possible to determine explicitly the eigenvalues and
the corresponding eigenfunctions of $A$ for a generic bounded domain $\Xi
\subset \erre^{2}$. However, for particular shapes of $\Xi $ the
eigensystem is well-known, e.g. for $\Xi $ being a rectangle. If $\Xi
=[0,L]\times \lbrack 0,H]$, one has 
$$
\Delta e_{m,n}=-\lambda _{m,n}e_{m,n}
$$
with 
$$
e_{m,n}(\xi_1,\xi_2)=
\cos\frac{m\pi \xi_1}{L} \sin\frac{n\pi\xi_2}{H}
$$
and 
$$
\lambda_{m,n}=\Big(\frac{m^2}{L^2} + \frac{n^2}{H^2}\Big)\pi^2.
$$
Moreover, $\lambda_{0,0}=0$ is an eigenvalue with unit eigenvector
$e_{0,0}=(LH)^{-1/2}$.
Hence $Ae_{m,n}=-(\lambda _{m,n}+\rho )e_{m,n}$.

\section{Analysis of the targeting formulation}
In the Hilbert space setting introduced in Section \ref{sec:reform},
assuming $B=I$ for simplicity, let us define the distance from the
target $y:[0,T]\to X$ as
$$
y(t,\xi )=h(\xi )-x(t,\xi ),
$$%
where $h\in L_2(\Xi )$ is the desired configuration of goodwill to reach
at time $T$. Then $y$ is the unique mild solution of the following
non-homogeneous linear Cauchy problem in $X$: 
\begin{equation}
\label{eq:tracking}
\left\{ 
\begin{array}{l}
\displaystyle y'=Ay(t)-u(t)+f, \\[6pt] 
y(0)=h-x_{0},%
\end{array}%
\right.  
\end{equation}%
with $f=-Ah$.

Problem (P2) can be rewritten as 
\begin{equation}
\inf_{u\in L_2([0,T];E)} \left[ \gamma
|y(T)|^{2}+\int_{0}^{T}|u(t)|^{2}\,dt \right],  \label{eq:P2-track}
\end{equation}%
subject to (\ref{eq:tracking}). Appealing again to the dynamic programming
principle, we can write the Riccati equation associated to problem (\ref%
{eq:P2-track}) as follows: 
\begin{equation}
\left\{ 
\begin{array}{l}
\displaystyle P'(t)=2AP(t)-P^{2}(t), \\[6pt] 
P(0)=\gamma I,%
\end{array}%
\right.  \label{eq:ricc21}
\end{equation}%
and its adjoint backward equation as 
\begin{equation}
\left\{ 
\begin{array}{l}
\displaystyle r'(t)+(A-P(t))r(t)+P(t)f=0, \\[6pt] 
r(T)=0.%
\end{array}%
\right.  \label{eq:ricc22}
\end{equation}%
Then one can uniquely solve the tracking problem in terms of the solution to
(\ref{eq:ricc21}) and (\ref{eq:ricc22}).

\begin{prop}
If the target function $h$ is such that $Ah\in L_2\left( \Xi
\right) $, then the optimal control problem (P2) admits a
unique optimal control $u_*$ given by the feedback law 
\[
u_*(t)=P\left( T-t\right) y_*(t)+r(t), 
\]%
where $y_*$ is the unique mild solution of the closed-loop
equation 
\begin{equation}
\label{eq:traj2}
\left\{ 
\begin{array}{l}
\displaystyle y'(t) = (A-P(T-t))y(t)-r(t)+f, \\ 
y(0) =h-x_{0},
\end{array}
\right.  
\end{equation}
and $P(\cdot)$, $r(\cdot)$ solve, respectively, equations
(\ref{eq:ricc21}) and (\ref{eq:ricc22}).  Moreover, the value function
$J^{\ast }$ is given by
\begin{eqnarray*}
J^* = J(u_*) &=& \cp{P(T)y(0)}{y(0)} + 2\cp{r(0)}{y(0)} \\
&&+\int_0^T \Big[ 2\cp{r(t)}{f} -|r(t)|^2\Big]\,dt.
\end{eqnarray*}
\end{prop}
\begin{proof}
  It follows essentially the same lines of the ``classical'' proof in
  \cite{DP-ICTP}, with the only difference that in this case $P_0$ is
  positive definite (hence we do not need conditions on the data of
  the problem to ensure well-posedness), and we have to take into
  account a non-homogeneous part in the state equation, and
  consequently of the adjoint equation (\ref{eq:ricc22}) accompanying
  the Riccati equation. The assumption $Ah\in L_2(\Xi)$ is
  equivalent to $f\in X$, so that a (unique) mild solution of
  (\ref{eq:tracking}) exists. For more details we refer to \cite{rep2}.
\end{proof}

Note that, in contrast to (P1), there always exists an optimal
solution for problem (P2), for any choice of the parameters and any initial
condition. Below we obtain an expression for the optimal trajectory and the
optimal control in terms of a basis of $L_2(\Xi)$, as we have done for
(P1). In particular, by projecting the Riccati equation (\ref{eq:ricc21}) on
the system $(e_{k})_{k\in \N}$, we obtain the infinite set of Cauchy
problems 
$$
\left\{ 
\begin{array}{l}
\displaystyle \frac{dp_{k}}{dt}=-2\lambda _{k}p_{k}(t)-p_{k}^{2}(t), \\[8pt] 
p_{k}(0)=\gamma ,%
\end{array}%
\right.
$$%
each of which admits the explicit solution 
$$
p_{k}(t)=-2\lambda _{k}+\frac{1}{(2\lambda _{k})^{-1}+C_{k}e^{-2\lambda
_{k}t}},
$$%
with $C_{k}=-\gamma (2\lambda _{k}(2\lambda _{k}+\gamma ))^{-1}$. As before,
we have set $p_{k}(\cdot ):=P(\cdot )e_{k}$.

The adjoint backward Cauchy problem (\ref{eq:ricc22}) can be solved
similarly: projecting on the system $(e_{k})_{k\in \N}$ we get 
$$
\left\{ 
\begin{array}{l}
\displaystyle \frac{dr_{k}}{dt}-\lambda
_{k}r_{k}(t)-p_{k}(T-t)r_{k}(t)+p_{k}(T-t)f_{k}=0, \\[8pt] 
r_{k}(T)=0,%
\end{array}%
\right.
$$%
where $r_{k}(\cdot ):=r(\cdot )e_{k}$ and $f_{k}:=fe_{k}$. Setting $\eta
_{k}(t):=r_{k}(T-t)$, one has 
$$
\left\{ 
\begin{array}{l}
\displaystyle \frac{d\eta _{k}}{dt}=\lambda _{k}\eta _{k}(t)-p_{k}(t)\eta
_{k}(t)+p_{k}(t)f_{k}, \\[8pt] 
\eta _{k}(0)=0,%
\end{array}%
\right.
$$%
These Cauchy problems can be solved explicitly, yielding 
\[
\eta _{k}(t)=\gamma f_{k}e^{\lambda
_{k}t+\int_{0}^{t}p_{k}(s)\,ds}+f_{k}\int_{0}^{t}e^{\lambda _{k}(t-\tau
)-\int_{\tau }^{t}p_{k}(s)\,ds}p_{k}(\tau ), 
\]%
and finally $r_{k}(t)=\eta _{k}(T-t)$.

The optimal trajectory can also be written explicitly by computing the
solution of the closed-loop equation. We again project the equation on the
system $(e_{k})_{k\in \N}$, obtaining 
$$
\left\{ 
\begin{array}{l}
\displaystyle \frac{dy^k}{dt}=(-\lambda_k-p_k(T-t))y^k(t)-r_k(t)+f_k, \\[8pt] 
y^k(0)=h_k-x_k(0),
\end{array}
\right.
$$
hence 
$$
y_*^k(t)=(h_{k}-x_{k}(0))e^{\int_{0}^{t}a_{k}(s)\,ds}+\int_{0}^{t}e^{\int_{\tau
}^{t}a_{k}(s)\,ds}(f_{k}-r_{k}(s))\,ds,
$$
where we have set $a_{k}(s):=-\lambda _{k}-p_{k}(T-s)$.
The optimal trajectory can now be written as 
$$
y_*(t,\xi)=\sum_{k=0}^{\infty }y^k_*(t)e_{k}(\xi ).
$$
Similarly, the optimal policy is given by 
$$
u_*(t,\xi)=\sum_{k=0}^\infty u^k_*(t)e_{k}(\xi)=
\sum_{k=0}^{\infty }(p_{k}(T-t)y^k_*(t)+r_k(t))e_k(\xi).
$$

\section{Discussion}
Our paper extends traditional marketing models of goodwill dynamics to
allow for spatially distributed advertising. Using an appropriate
Hilbert space reformulation, we map a problem of profit maximization
for a monopoly firm into an optimal control problem in infinite
dimensions and discuss existence and characterization of its optimal
solutions. In some simple, but still realistic situations, the optimal
strategy can be obtained in closed form. Our analysis provides a
tractable model to investigate how advertising funds should be
distributed over multiple markets.

Some important issues to be considered as follow-up investigations
are, for instance, an assessment of the empirical adequacy of our
model, as well as the analysis of dynamic interaction between firms
advertising in multiple markets. Such analysis would present an
important complement to the existing literature on single-market
advertising competition.

\section*{Acknowledgments}
The first author gratefully acknowledges financial support by the
DAAD, the LuRo foundation through grant MAC-PMB-03, the SFB 611
(Bonn), and the Max-Planck-Institut f\"ur Mathematik (Leipzig) through
an IPDE fellowship. Part of this work was carried out during a visit
of the first author at the Department of Mathematics of the University
of Trento, whose hospitality is gratefully acknowledged.

\def\cprime{$'$}

\end{document}